\documentclass[12pt,twoside]{article}
\usepackage{amsmath,amssymb,mathrsfs,graphicx,hyperref}

\newtheorem{theorem}{Theorem}[section]

\newtheorem{corollary}[theorem]{Corollary}

\newtheorem{rmrk}[theorem]{Remark}

\DeclareMathAlphabet{\mathbfit}{OML}{cmm}{b}{it}

\makeatletter
\@addtoreset{equation}{section}
\makeatother

\newenvironment{remark}
{\begin{rmrk} \em}
{\end{rmrk}}

\newcommand{\R} {\mathbb{R}}

\newcommand{\Z} {\mathbb{Z}}
\newcommand{\N} {\mathbb{N}}

\newcommand{\ds} {\displaystyle}

\newcommand{\into} {\longrightarrow}

\renewcommand{\a} {\alpha}
\renewcommand{\b} {\beta}
\renewcommand{\o} {\omega}
\renewcommand{\P} {\mathbb{P}}
\renewcommand{\tilde} {\widetilde}
\newcommand{\mq} {\overline{m}_q}
\newcommand{\bq} {p}

\newcommand{\tinyspace}{\hspace{0.1em}}

\newcommand{\cadlag}{\ifmmode\mathcal D\else c\`adl\`ag \fi}

\newcommand{\id}{\mathrm{id}}

\newcommand{\dconv}{\xrightarrow{\:\: \mathrm{d} \:\: }}   
\newcommand{\asconv}{\xrightarrow{\:\: \mathrm{a.s.} \:\: }}   

\newcommand{\E}[1]{\mathbb E[#1]}   
\newcommand{\inDJ}[1]{
    \IfEqCase{#1}{
        {1}{\qquad \text{in } (\cadlag,J_1)}
        {2}{\qquad \text{in } (\cadlag,J_2)}
    }[\PackageError{inDJ}{Undefined option to inDJ: #1}{}]%
}

\newcommand{\inDM}[1]{
    \IfEqCase{#1}{
        {1}{\qquad \text{in } (\cadlag,M_1)}
        {2}{\qquad \text{in } (\cadlag,M_2)}
    }[\PackageError{inDM}{Undefined option to inDM: #1}{}]%
}
\newcommand{\inDPJ}[1]{
    \IfEqCase{#1}{
        {1}{\qquad \text{in } (\cadlag^{+},J_1)}
        {2}{\qquad \text{in } (\cadlag^{+},J_2)}
    }[\PackageError{inDPJ}{Undefined option to inDPJ: #1}{}]%
}
\newcommand{\inDPM}[1]{
    \IfEqCase{#1}{
        {1}{\qquad \text{in } (\cadlag^{+},M_1)}
        {2}{\qquad \text{in } (\cadlag^{+},M_2)}
    }[\PackageError{inDPM}{Undefined option to inDPM: #1}{}]%
}

\newcommand{\RW}[1][n]{S_{#1}}
\newcommand{\fluidscaledRW}[1][n]{\bar{S}^{(#1)}}
\newcommand{\diffusionscaledRWnomean}[1][n]{\hat{S}^{(#1)}}
\newcommand{\diffusionscaledRW}[1][n]{\tilde{S}^{(#1)}}
\newcommand{\diffusionlimitRW}[1][\beta]{W^{(#1)}}  

\newcommand{\environment}[1][n]{\omega_{#1}}
\newcommand{\fluidscaledenvironment}[1][n]{\bar{\omega}^{(#1)}}
\newcommand{\diffusionscaledenvironmentnomean}[1][n]{\hat{\omega}^{(#1)}}
\newcommand{\diffusionscaledenvironment}[1][n]{\tilde\omega^{(#1)}}
\newcommand{\diffusionlimitenvironment}[1][\alpha]{Z^{(#1)}}  

\newcommand{\flight}[1][n]{Y_{#1}}
\newcommand{\fluidscaledflight}[1][n]{\bar{Y}^{(#1)}}
\newcommand{\diffusionscaledflight}[1][n]{\tilde{Y}^{(#1)}}
\newcommand{\diffusionscaledflightnomean}[1][n]{\hat{Y}^{(#1)}}

\begin{document}

\title{\textbf{Discrete- and continuous-time random walks in 1D L\'evy 
random medium}}

\author{
\scshape
Marco Lenci\thanks{
Dipartimento di Matematica, Universit\`a di Bologna,
Piazza di Porta San Donato 5, 40126 Bologna, Italy; and
Istituto Nazionale di Fisica Nucleare, Sezione di Bologna, Via Irnerio 46,
40126 Bologna, Italy; E-mail:\href{mailto:marco.lenci@unibo.it}{
\texttt{marco.lenci@unibo.it}}.}
\vspace*{7pt}
}

\date{In memory of Francesca Romana Nardi
 \\[5pt]
January 2022}

\maketitle

\begin{abstract}
A L\'evy random medium, in a given space, is a random point process 
where the distances between points, a.k.a.\ \emph{targets}, are long-tailed. 
Random walks visiting the targets of a L\'evy random medium have been used to 
model many (physical, ecological, social) phenomena that exhibit superdiffusion 
as the result of interactions between an agent and a sparse, complex environment. 
In this note we consider the simplest non-trivial L\'evy random medium, a 
sequence of points in the real line with i.i.d.\ long-tailed distances between
consecutive targets. A popular example of a continuous-time random walk in 
this medium is the so-called \emph{L\'evy-Lorentz gas}. We give an account 
of a number of recent theorems on generalizations and variations of such 
model, in discrete and continuous time.

  
\end{abstract}

\section{Introduction}
\label{sec-intro}

In this note we give an overview of recent rigorous results on random walks 
(RWs) in random medium on the real line. The random medium is given by 
a point process $\o = (\o_k, \, k \in \Z) \subset \R$, where 
$\o_0=0$ and the distances $\zeta_k := \o_k - \o_{k-1}$ 
between consecutive points are positive i.i.d.\ random variables with a long 
tail. By this we mean that the variance of $\zeta_k$ is infinite. The points
$\o_k$ will be henceforth called \emph{targets}.
For reason that will be better clarified below we refer to $\o$ as
a \emph{L\'evy random medium}. 

We consider two types of RWs on $\R$ related to $\o$.
To define them we introduce the 
auxiliary process $S = (S_n, \, n \in \N)$, a $\Z$-valued RW with $S_0=0$ 
and independent increments. We postulate that $S$ is independent of 
$\o$ and call it the \emph{underlying random walk}. The first process of 
interest is $Y \equiv Y^\o := (Y_n, \, n \in \N)$, where $Y_n := \o_{S_n}$. 
This is the discrete-time RW (DTRW) that 
``jumps'' on the targets of $\o$ as determined by $S$. For example, 
if $S$ produces the realization $(0, 2, -3, \ldots)$, the walker $Y$
starts at the origin, then jumps to the second target to the right of 0, then to 
the third target to the left of 0, etc. 
The second process of interest is $X \equiv X^\o := ( X(t), \, t \ge 0)$, 
the continuous-time RW 
(CTRW) defined as the unit-speed interpolation of $Y$. This means that
the walker $X$ visits all the points $Y_n$, ordered by $n$, but ``walking'' 
with unit speed rather than jumping. For instance, for $S$ as in the 
above example, $X$ starts at the origin and moves with velocity $1$ 
until it reaches $Y_1 = \o_2$, then it instantaneously turns its velocity 
to $-1$ and moves until it reaches $Y_2 = \o_{-3}$, and so on. 

In the case where the underlying RW $S$ is simple and symmetric, the
process $X$ is generally referred to as the \emph{L\'evy-Lorentz gas}, 
after Barkai, Fleurov and Klafter introduced it in the physical literature in 2000 
\cite{bfk}. The L\'evy-Lorentz gas has been used since as a simple model 
for a number of phenomena exhibiting superdiffusion, i.e., diffusion at a
faster speed than square root of time. They include transport in porous 
media, disordered optical media (such as \emph{L\'evy glasses} \cite{bbw}), 
nanowires, etc.; see \cite{bfk, abg, bcv, vbb} and references therein.

In the physical literature, DTRWs, respectively CTRWs, whose distributions 
of jumps, respectively inertial stretches, are long-tailed, 
are often called \emph{L\'evy flights}, respectively \emph{L\'evy 
walks}. L\'evy flights and walks, in regular or random
media, have been employed as models for anomalous diffusion in a 
wide range of situations, from the physical to the biological and social
sciences \cite{krs, zdk}. A rigorous mathematical treatment of these systems 
has only been given for the simplest of them, mostly on regular media. In 
real-world applications, however, the anomalous behavior of a certain 
diffusing quantity is seldom due to a special law governing the diffusing agent 
\emph{per se}, but rather to the interaction between the agent and an irregular 
medium (e.g., a photon in a L\'evy glass, a signal in a small-world network, an 
animal foraging where food is scarce, etc.). Hence the interest in \emph{L\'evy 
media}, namely, media that induce superdiffusive behavior, such as the 
random point process $\o$ defined earlier.

A fair amount of mathematical work on the processes $X$ and $Y$ has 
been done in recent years by different authors \cite{bcll, ms, blp, sbblm, z}. 
In particular, in a number of cases, limit theorems have been proved 
for suitable rescalings of either process. In some instances, the 
convergence of moments has been proved as well. The purpose of this 
note is to present these results 
in a concise, unified manner. For reasons of space and self-consistency, 
we will neglect interesting work 
by Artuso and collaborators on yet another type of RW related to the 
L\'evy-Lorentz gas, a persistent RW in an \emph{averaged} medium
\cite{acor, roac}.

In Section \ref{sec-dtrw} we present results on the DTRW $Y$ and in
Section \ref{sec-ctrw} on the CTRW $X$. Understandably, the results 
depend on the assumptions on the random medium $\o$ and the
underlying RW $S$. Major differences occur depending on whether 
$\zeta_k$, the distance between two consecutive targets, has infinite
variance but finite mean, or infinite mean (and thus infinite variance),
so we consider these cases in different subsections. No proofs are
given, but references are placed throughout.

\subsection{General notation}

Throughout the paper we denote by $\P$ the probability law that governs 
the whole system, both the random medium $\o$ and the underlying 
RW $S$ (that latter playing the role of the random dynamics, as it
``drives'' both $X$ and $Y$, in a given $\o$). $\P$ is called the 
\emph{annealed law} and we denote by $\mathbb{E}$ its expectation. 
For a fixed $\o$, the 
conditional probability $P_\o := \P( \cdot | \o)$ is called the 
\emph{quenched law} relative to the realization $\o$ of the medium. 
We denote by $E_\o$ its expectation. A limit theorem, such as the CLT
or the Invariance Principle, etc., relative to $\P$ is referred to as an 
\emph{annealed limit theorem}. One speaks instead of a \emph{quenched 
limit theorem} if the result is proved w.r.t.\ $P_\o$, for $\P$-a.e.\
$\o$.%
\footnote{Since we have neither introduced the measurable space 
$(\Omega, \mathcal{A})$ where $\P$ is defined, nor declared that 
$\o$ are elements of $\Omega$, mathematical formality requires that 
we define the phrase ``$\P$-almost every $\o$''. The counterimages 
(equivalently, level sets) of the process $\o$ form a partition of
$\Omega$. We assume this partition to be \emph{measurable} in the 
sense of Rohlin \cite{r}. Now, a property is said to hold for $\P$-a.e.\ $\o$ 
if the values of $\o$ which do not satisfy the property correspond to
elements of the partition whose union has zero $\P$-measure. 
Incidentally, the existence of such a measurable partition is what 
guarantees that $\P( \cdot | \o)$ is well-defined (for $\P$-a.e.\ $\o$).}
In what follows, for the most part, we present annealed and quenched 
results in different subsections.

We indicate with $\xi_n := S_n - S_{n-1}$ ($n \in \Z^+$) the i.i.d.\ increments
of the underlying RW, whose drift is denoted $\nu := \E{\xi_n}$, if it exists
in $\R \cup \{ \pm \infty \}$. The mean distance between the targets of $\o$
is denoted $\mu := \E{\zeta_k}$. Since $\zeta_k>0$, $\mu$ always exists in 
$\R \cup \{ +\infty \}$. Let us recall that all $\xi_n$ ($n \in \Z^+$) and $\zeta_k$ 
($k \in \Z$) are independent.

\section{Discrete-time random walk}
\label{sec-dtrw}

In this section we consider the asymptotic behavior of the DTRW 
$Y \equiv Y^\o$, under a number of different assumptions on the distributions of 
$\zeta_1$ and $\xi_1$. 

\subsection{Finite mean distance between targets, quenched theorems}
\label{subs-y-fin-q}

We start with the results of \cite{bcll} on the quenched version of $Y$, 
which only require a very simple condition on the medium, 
$\mu = \E{\zeta_1} < \infty$, that is, the 
mean distance between neighboring targets is finite. The 
assumptions on the underlying RW are instead as follows:
\begin{itemize}
\item the increment $\xi_1$ of $S$ is symmetric, i.e., $\P(\xi_1=j) = 
  \P(\xi_1=-j)$, for all $j \in \N$;
\item its distribution is unimodal, i.e., $j \mapsto \P(\xi_1=j)$ is non-increasing 
  for $j \in \N$;
\item it has finite variance: $V_\xi := \E{\xi_1^2} < \infty$.  
\end{itemize}
The authors prove a quenched CLT for $Y$ \cite[Thm.~1]{bcll}:
\begin{theorem} \label{thm-1}
  Assume the above conditions, most notably $\mu < \infty$. Then, as $n \to \infty$,
  \begin{displaymath}
    \frac{Y_n}{\sqrt{n}} \dconv \mathcal{N} \big(0, \mu^2 \, V_\xi \big) ,
  \end{displaymath}
  w.r.t.\ $P_\o$, for $\P$-a.e.\ $\o$. Here $\mathcal{N}(0, \mu^2 \, V_\xi)$ is a Gaussian 
  variable with mean 0 and variance $\mu^2 V_\xi$.
\end{theorem}
Obviously, a quenched distributional limit theorem with the same limit for a.e.\
quenched law implies the annealed version of the same theorem:
\begin{corollary} \label{cor-thm-1}
  The limit in the statement of Theorem \ref{thm-1} holds w.r.t.~$\P$ as well.
\end{corollary}

Convergence is known for the quenched moments of $Y_n/\sqrt{n}$ as well, at 
least of lower order. Let $\bq := \sup \{ q\ge 0 \:|\: \E{|\xi_1|^q} < \infty \}$. 
By the assumption on $V_\xi$, $\bq \ge 2$. For all $q \in \R^+$, 
denote by 
\begin{equation} \label{mq}
  \mq := \sqrt{ \frac{2^q} \pi } \, \Gamma \! \left( \frac{q+1} 2 \right)
\end{equation}
the $q$-th absolute moment of the standard Gaussian $\mathcal{N}(0,1)$ 
(here $\Gamma$ is the usual Gamma function). It is not hard to show that,
at least for all $q < \bq$,
\begin{equation}
  \lim_{n \to \infty} \frac{\E{S_n}}{n^{q/2}} = V_\xi^{q/2} \, \mq .
\end{equation}
The following is a reformulation of Theorem 2 of \cite{bcll}.
\begin{theorem} \label{thm-2}
  Under the above assumptions and notation, fix $q \in (0, \bq)$. For
  a.a.\ $\o$, 
  \begin{displaymath}
    \lim_{n \to \infty} \frac{E_\o[Y_n]}{n^{q/2}} = \mu^q \, V_\xi^{q/2} \, \mq .
  \end{displaymath}
\end{theorem}
Observe that $\mu^q \, V_\xi^{q/2} \, \mq$ is 
the $q$-th absolute moment of $\mathcal{N}(0, \mu^2 \, V_\xi)$, so Theorem
\ref{thm-2} is consistent with Theorem \ref{thm-1}.

\subsection{Finite mean distance between targets, annealed theorems}
\label{subs-y-fin-a}

In the next two subsections we report the functional limit theorems of 
\cite{sbblm}. We refer the reader to \cite{w} for background material on 
stable laws, L\'evy processes, Skorokhod topologies, etc. All distributional 
convergences in these subsections are meant w.r.t.\ $\P$, that is, we are 
considering annealed functional limit theorems.

We assume that $\zeta_1$ is in the normal basin of attraction of a 
$\a$-stable distribution, with $\a \in (1,2]$. This means that 
$\mu = \E{\zeta_1} <\infty$ and 
\begin{equation} \label{tilde-z1}
  \frac1 {n^{1/\a}} \sum_{i=1}^n (\zeta_i-\mu) \dconv \tilde{Z}^{(\a)}_1, 
\end{equation}
as $n \to \infty$, for some $\a$-stable variable $\tilde{Z}^{(\a)}_1$ (whose 
skewness index%
\footnote{This is the parameter that, in virtually all textbooks on stable 
variables (such as \cite{w}) is denoted $\beta \in [-1,1]$. In this paper 
$\beta$ is used for the stability index of $\xi_1$.}
 must then be 0). As for the underlying RW, we assume 
$\xi_1$ is in the normal basin of attraction of an $\b$-stable distribution, 
with $\b \in (0,1) \cup (1,2]$. We must distinguish two cases, depending 
on whether $\nu = \E{\xi_1}$ exists and differs from 0, or otherwise.
\begin{itemize}

\item If $\b \in (0,1)$, or $\b \in (1,2]$ and $\nu=0$, we assume that there exists 
a $\b$-stable variable $W^{(\b)}_1$ such that, as $n \to \infty$, 
\begin{equation} \label{w1}
  \frac1 {n^{1/\b}} \sum_{i=1}^n \xi_i \dconv W^{(\b)}_1.
\end{equation}

\item If $\b \in (1,2]$ and $\nu \ne 0$, we assume that there exists 
a $\b$-stable variable $\tilde{W}^{(\b)}_1$ such that
\begin{equation} \label{tilde-w1}
  \frac1 {n^{1/\b}} \sum_{i=1}^n (\xi_i-\nu) \dconv \tilde{W}^{(\b)}_1.
\end{equation}

\end{itemize}

To state the results of this section, we need two spaces of functions with jump 
discontinuities. In what follows, we denote by $\cadlag^+$ the space of 
c\`adl\`ag%
\footnote{I.e., right-continuous with left limits existing everywhere. C\`agl\`ad
means left-continuous with right limits everywhere.}
 functions $\R^+ \into \R$ and by $\cadlag$ the space of functions
$\R \into \R$ whose restriction to $[0, +\infty)$, respectively $(-\infty,0]$, is
c\`adl\`ag, respectively c\`agl\`ad. 

Let $(\tilde{Z}^{(\a)}_{\pm}(s), \, s \ge 0)$ be two i.i.d.\ \cadlag $\a$-stable 
L\'evy processes such that $\tilde{Z}^{(\a)}_{\pm}(0)=0$ and 
$\tilde{Z}^{(\a)}_{\pm}(1)$ is distributed like $\tilde{Z}^{(\a)}_1$, introduced in
(\ref{tilde-z1}) (these conditions uniquely determine the distribution of the 
processes), and set
\begin{equation} \label{tilde-z}
  \tilde{Z}^{(\a)}(s) := \left\{
  \begin{array}{lll}
    Z^{(\a)}_{+}(s) \,, && s \ge 0; \\
    -Z^{(\a)}_{-}(-s) \,, && s < 0.
  \end{array}
  \right.
\end{equation}
By construction, every realization $\tilde{Z}^{(\a)}$ belongs to $\cadlag$, and
so do the realizations
\begin{equation}
  \fluidscaledenvironment (s) := \frac1n \left\{
  \begin{array}{lll}
    \environment[\lfloor n s \rfloor] \,, && s \ge 0; \\
    \environment[\lceil n s \rceil] \,, && s < 0,
  \end{array}
  \right.
\end{equation}
\begin{equation} 
  \diffusionscaledenvironment (s) := \frac1 {n^{1/\a}} \left\{
  \begin{array}{lll}
    \sum_{i=1}^{\lfloor ns \rfloor} \, (\zeta_i - \mu) \,, && s \ge 0; \\[3pt]
    -\sum_{i= \lceil (n-1) s \rceil}^0 \, (\zeta_i - \mu) \,, && s < 0,
  \end{array}
  \right.
\end{equation}
defining the processes $( \fluidscaledenvironment (s), \, s\in \R)$ and 
$( \diffusionscaledenvironment (s), \, s\in \R)$. The single-variable 
convergence (\ref{tilde-z1}) entails functional convergence of these 
processes:  as $n\to \infty$, $\fluidscaledenvironment \asconv \mu 
\tinyspace\id$ and $\diffusionscaledenvironment \dconv \tilde{Z}^{(\a)}$, 
relative to the Skorokhod topology $J_1$. From now on, we will write `in
$(\cadlag, J_1)$' for short.

We now introduce continuous-argument processes for the dynamics. 
\begin{itemize}

\item In the case $\b \in (0,1)$, or $\b \in (1,2]$ and $\nu=0$, we denote
by $(W^{(\b)}(t), \, t \ge 0)$ a \cadlag $\b$-stable L\'evy process whose 
distribution is uniquely determined by the conditions that $W^{(\b)}(0)=0$ 
and $W^{(\b)}(1)$ be distributed like $W^{(\b)}_1$; cf.\ (\ref{w1}). Also
define $( \diffusionscaledRWnomean (t) , \, t \ge 0)$ via
\begin{equation} \label{hat-S}
\diffusionscaledRWnomean (t) := \frac{ \RW[\lfloor n t \rfloor]}{n^{1/\b}}.
\end{equation}
It follows from (\ref{w1}) that $\diffusionscaledRWnomean \dconv 
\diffusionlimitRW$, in $(\cadlag^+, J_1)$, as $n \to \infty$.

\item In the case $\b \in (1,2]$ and $\nu \ne 0$, we consider 
$\tilde{W}^{(\b)}$, defined exactly as $W^{(\b)}$ above but 
with $\tilde{W}^{(\b)}_1$ in place of $W^{(\b)}_1$; cf.\ (\ref{tilde-w1}).
In lieu of (\ref{hat-S}) we define two processes:
\begin{equation} 
  \fluidscaledRW (t) := \frac{ \RW[\lfloor n  t \rfloor]}n \,, \qquad
  \diffusionscaledRW(t) := \frac{\sum_{i=1}^{\lfloor nt \rfloor} \,
  (\xi_{i} - \nu)} {n^{1/\b}} .
\end{equation}
All these processes take values in $\cadlag^+$. By 
(\ref{tilde-w1}), $\fluidscaledRW \asconv \nu \tinyspace\id$
and $\diffusionscaledRW \dconv \widetilde{W}^{(\b)}$, in  
$(\cadlag^+, J_1)$, as $n \to \infty$.

\end{itemize}
The following result extends Theorem 2.3 of \cite{sbblm}.
\begin{theorem} \label{thm-3}
  Under the above assumptions, in particular $\a \in (1,2]$, the following
  convergences hold, w.r.t.~$\P$:
  \begin{itemize}
    \item[(a)] If $\b \in (0,1)$, or $\b \in (1,2]$ with $\nu=0$, 
    let $\diffusionscaledflightnomean(t) := \ds \frac{ \flight[\lfloor nt \rfloor] }
    { n^{1/\b} }$, for $t\ge 0$. As $n\to\infty$,
    \begin{displaymath}
      \diffusionscaledflightnomean \dconv \mu \tinyspace \diffusionlimitRW
      \qquad \text{in } (\cadlag^+, J_1).
    \end{displaymath}
    
    \item[(b)] If $\b \in (1,2]$ with $\nu \ne 0$, let
    $\fluidscaledflight(t) :=  \ds \frac{ \flight[\lfloor nt \rfloor] }{ n }$,
    for $t\ge 0$. As $n\to\infty$,
    \begin{displaymath}
     \fluidscaledflight \dconv \mu \nu \tinyspace \id
      \qquad \text{in } (\cadlag^+, J_1).
    \end{displaymath}
  \end{itemize}
\end{theorem}

\begin{remark}
  The statement of \cite[Thm.~2.3]{sbblm} does not include the cases
  $\a=2$ and/or $\b = 2$, because the authors were mostly interested
  in \emph{bona fide} L\'evy media and L\'evy flights in them. The proof of the 
  theorem, however, works verbatim if the assumptions are generalized to 
  include $\zeta_1$ and/or $\xi_1$ in the normal domain of attraction of 
  a 2-stable distribution, i.e., a Gaussian. In this case, of course, 
  $\tilde{Z}^{(2)}_{\pm}$, $W^{(2)}_1$ and $\tilde{W}^{(2)}_1$ are
  Brownian motions. The same remark holds for Theorems \ref{thm-4} and
  \ref{thm-5} below.
 \end{remark} 
 
\begin{remark}
  Theorem \ref{thm-3} mostly supersedes Theorem 2.1 of 
  \cite{ms} (which is stated for the case where $S$ is simple and 
  symmetric), but not quite, since the hypothesis on $\zeta_1$ there is that
  $\P(\zeta_1 > x) \approx x^{-\bq}$, for $x \to +\infty$, with $\bq \ge 1$.
  This is weaker than asking that $\zeta_1$ be in the \emph{normal} 
  domain of attraction of an $\a$-stable distribution, with $\a \in (1,2]$.
  For example it includes the case $\P(\zeta_1 > x) = C x^{-2}$, where 
  $\zeta_1$ is in the domain, but not normal domain of attraction of a
  Gaussian \cite[Ex.~5.10]{j}. Also, the assertion of Theorem \ref{thm-3} 
  is of course much stronger than that of Corollary \ref{cor-thm-1}, 
  but the hypothesis on the medium for the latter is much weaker: simply
  $\mu < \infty$.
\end{remark}

The case \emph{(b)} of the above theorem is the case where $Y$ has a drift. 
Understandably,
the scaling rate of $\flight[\lfloor n\cdot \rfloor]$ is $n$ (one says that the 
process is \emph{ballistic}) 
and the convergence is to a deterministic function. It is therefore natural to 
study the fluctuations around the deterministic limit.
\begin{theorem} \label{thm-4}
  Under the same assumptions and notation as Theorem \ref{thm-3}, if 
  $\a, \b \in (1,2]$, the following convergences hold, w.r.t.~$\P$:
  \begin{itemize}
    \item[(a)] If $\b < \a$, let $\diffusionscaledflight(t) := \ds \frac{ 
    n(\fluidscaledflight(t) - \mu \nu \tinyspace t)} { n^{1/\b} } = \frac{ 
    \flight[\lfloor nt \rfloor] - n \mu \nu \tinyspace t } { n^{1/\b} }$. As $n\to\infty$,
    \begin{displaymath}
      \diffusionscaledflight \dconv \mu \tilde{W}^{(\b)}
      \qquad \text{in } (\cadlag^+, J_1).
    \end{displaymath}
 
     \item[(b)] If $\b > \a$, let $\diffusionscaledflight(t) := \ds \frac{ 
    n(\fluidscaledflight(t) - \mu \nu \tinyspace t)} { n^{1/\a} } = \frac{ 
    \flight[\lfloor nt \rfloor]  - n \mu \nu \tinyspace t } { n^{1/\a} }$. As $n\to\infty$,
    \begin{displaymath}
      \diffusionscaledflight \dconv \mathrm{sgn}(\nu) \, |\nu|^{1/\a} \, 
      \tilde{Z}^{(\a)}_+ 
      \qquad \text{in } (\cadlag^+, J_2).
    \end{displaymath}
    
    \item[(c)] If $\b = \a$, let $\diffusionscaledflight(t) := \ds \frac{ 
    n(\fluidscaledflight(t) - \mu \nu \tinyspace t)} { n^{1/\b} } = \frac{ 
    \flight[\lfloor nt \rfloor]  - n \mu \nu \tinyspace t } { n^{1/\b} }$. As $n\to\infty$,
    \begin{displaymath}
      \diffusionscaledflight \dconv \mu \tilde{W}^{(\b)} +
      \mathrm{sgn}(\nu) \, |\nu|^{1/\a} \, \tilde{Z}^{(\a)}_+ 
      \qquad \text{in } (\cadlag^+, J_2).
    \end{displaymath}
    where $\tilde{W}^{(\b)}$ and $\tilde{Z}^{(\a)}_+$ are two independent
    processes, defined as in (\ref{tilde-w1}) ff.
  \end{itemize}  
\end{theorem}

\begin{remark} \label{rk-a}
  The limits in \emph{(b)} and \emph{(c)} involve the unusual
  Skorokhod topology $J_2$ \cite[\S11.5]{w}. But this is the strongest
  amongst the classical Skorokhod topologies relative to which such
  limits hold, cf. Remark 2.11 of \cite{sbblm}. However, see Remark A.2
  in the same paper. 
\end{remark}

\subsection{Infinite mean distance between targets, annealed theorems}
\label{subs-y-inf-a}

In this subsection we assume that $\zeta_1$ is in the normal basin of attraction 
of a $\a$-stable distribution with $\a \in (0,1)$. Since $\E{\zeta_1} = \infty$,
this means that, for $n \to \infty$,
\begin{equation} \label{z1}
  \frac1 {n^{1/\a}} \sum_{i=1}^n \zeta_i \dconv Z^{(\a)}_1, 
\end{equation}
for some $\a$-stable variable $Z^{(\a)}_1$ (whose skewness index is 1,
since $\zeta_i > 0$). Out of $Z^{(\a)}_1$, we construct continuous-argument 
processes $(Z^{(\a)}_{\pm}(s), \, s \ge 0)$ and $(Z^{(\a)}(s), \, s \in \R)$ in
complete analogy with the previous case; cf.\ (\ref{tilde-z}). $Z^{(\a)}$ 
takes values in $\cadlag$, and the same is true for
\begin{equation} \label{hat-o-n}
  \diffusionscaledenvironmentnomean (s) := \frac1 {n^{1/\a}} \left\{
  \begin{array}{lll}
    \environment[\lfloor n s \rfloor] \,, && s \ge 0, \\
    \environment[\lceil n s \rceil] \,, && s < 0.
  \end{array}
  \right.
\end{equation}
It is a basic fact that, as $n \to \infty$, $\diffusionscaledenvironmentnomean 
\dconv \diffusionlimitenvironment$, in $(\cadlag, J_1)$. 

As the underlying RW, we maintain the same assumptions and notation
as in \S\ref{subs-y-fin-a}, recalling that the fundamental assumption is that
$\xi_1$ is in the normal basin of attraction of an $\b$-stable distribution, 
with $\b \in (0,1) \cup (1,2]$. Once again, $\nu$ denotes the expectation of
$\xi_1$, when defined.

The following theorem comprises and extends Theorems 2.1 and 2.2 of
\cite{sbblm}:
\begin{theorem} \label{thm-5}
  Under the above assumptions, in particular $\a \in (0,1)$, the following
  convergences hold, w.r.t.~$\P$:
  \begin{itemize}
    \item[(a)] If $\b \in (0,1)$, or $\b \in (1,2]$ with $\nu=0$, 
    let $\diffusionscaledflightnomean(t) := \ds \frac{ \flight[\lfloor nt \rfloor] }
    { n^{1/\a\b} }$, for $t\ge 0$. As $n\to\infty$, the finite-dimensional distributions 
    of $\diffusionscaledflightnomean$ converge to those of 
    $\diffusionlimitenvironment \circ \diffusionlimitRW$. This means that, for all 
    $m \in \Z^+$ and $t_1,\ldots,t_m \in \R^+$,
    \begin{displaymath} 
      \big( \diffusionscaledflightnomean (t_1) , \ldots , 
      \diffusionscaledflightnomean (t_m) \big) \dconv
      \big( \diffusionlimitenvironment(\diffusionlimitRW (t_1)) , \ldots, 
      \diffusionlimitenvironment(\diffusionlimitRW (t_m)) \big).
    \end{displaymath} 

    \item[(b)] If $\b \in (1,2]$ with $\nu \ne 0$, let 
    $\diffusionscaledflightnomean(t) := \ds \frac{ \flight[\lfloor nt \rfloor] }
    { n^{1/\a} }$, for $t\ge 0$. As $n\to\infty$, 
    \begin{displaymath}
     \diffusionscaledflightnomean \dconv \mathrm{sgn}(\nu) \, |\nu|^{1/\a}
      \, Z^{(\a)}_+ 
      \qquad \text{in } (\cadlag^+, J_2).
    \end{displaymath}
  \end{itemize}
\end{theorem}

\begin{remark} \label{rk-d}
  The convergence of the finite-dimensional distributions in assertion 
  \emph{(a)} is certainly a weak form of convergence, but it is morally the
  best one can do, given that $\diffusionlimitenvironment \circ \diffusionlimitRW$
  is not \cadlag with positive probability; see the comments after Theorem 2.2 
  of \cite{blp}. As for assertion \emph{(b)}, the considerations of Remark
  \ref{rk-a} apply here too.
\end{remark}

\section{Continuous-time random walk}
\label{sec-ctrw}

In this section we deal with the CTRW $X \equiv X^\o$, again under various 
assumptions, depending on the papers we report on.

\subsection{Finite mean distance between targets, quenched theorems}

Once more, we start by presenting a result by \cite{bcll}, namely the 
quenched CLT for $X$ \cite[Thm.~1]{bcll}. The assumptions are the same 
as in \S\ref{subs-y-fin-q} above: the mean distance $\mu$ between targets 
is finite and the underlying RW has symmetric, unimodal, finite-variance
increments $\xi_n$. Let us recall in particular the notation $V_\xi := \E{\xi_1^2}$. 
\begin{theorem} \label{thm-6}
  Under the above assumptions, most notably $\mu < \infty$, let $M_\xi := 
  \E{|\xi_1|}$ denote the first absolute moment of the underlying RW. Then, as 
  $t \to \infty$ and w.r.t.\ $P_\o$, for a.a.\ $\o$, 
  \begin{displaymath}
    \frac{X(t)}{\sqrt{t}} \dconv \mathcal{N} \! \left(0, \mu \frac{V_\xi} {M_\xi} \right) .
  \end{displaymath}
  Here, once again, $\mathcal{N}(0, \cdot\, )$ is a centerd Gaussian variable with 
  the specified variance.
\end{theorem}
The annealed CLT follows immediately:
\begin{corollary} \label{cor-thm-6}
  The limit in the statement of Theorem \ref{thm-6} holds w.r.t.~$\P$ as well.
\end{corollary}

A recent preprint of Zamparo \cite{z} claims the convergence of \emph{all}
quenched moments of $X(t)/\sqrt{t}$, under the additional assumption that the 
underlying RW is simple and symmetric (implying that $X$ is the \emph{bona
fide} L\'evy-Lorentz gas). Recall the notation $\mq$ for 
the $q$-th absolute moment of the standard Gaussian, cf.\ (\ref{mq}).
\begin{theorem} \label{thm-7}
  Assume that $\mu < \infty$ and $S$ is a simple symmetric RW. Then,
  for a.a.\ $\o$, 
  \begin{displaymath}
    \lim_{t \to \infty} \frac{E_\o[ |X(t)|^q ]} {t^{q/2}} = \mu^{q/2} \, \mq .
  \end{displaymath}
\end{theorem}

\begin{remark} \label{rk-b}
  When $S$ is simple and symmetic, $V_\xi = M_\xi = 1$, so Theorem 
  \ref{thm-7} is consistent with Theorem \ref{thm-6}, showing that $X$ is
  completely diffusive in this case.
\end{remark}

Theorem \ref{thm-7} descends from another result of independent interest, 
concerning the \emph{large deviations} of $X$, namely, events of the type 
$\{ |X(t)| > at \}$, for $a>0$. Since $X(t)$ is centered and scales like $\sqrt{t}$, 
the probability of such ``ballistic events'' is expected to be exceedingly small.
In \cite[Thm.~2.3]{z} it is proved that this probability vanishes like a stretched 
exponential. We report such result here:
\begin{theorem} \label{thm-8}
  Under the same assumptions as in Theorem \ref{thm-7}, there exists 
  $\kappa > 0$ such that, for all $a \in (0,1]$, the limit
  \begin{displaymath}
    \limsup_{t \to \infty} \frac1 {\sqrt{at}} \log P_\o ( |X(t)| > at ) \le -\kappa
  \end{displaymath}
  holds for a.e.~$\o$.
\end{theorem}

\subsection{Finite mean distance between targets, annealed theorems}

Apart from recalling Corollary \ref{cor-thm-6}, which establishes the
annealed CLT for $X$ under the assumptions of \cite{bcll} ($\mu < \infty$
and $S$ has symmetric, unimodal, finite-variance increments, cf.\ \S2.1),
in this section we present the results of \cite{z} on the moments and 
large deviations of the annealed version of $X$. 

The assumptions for this part are stronger than for Theorems \ref{thm-7}
and \ref{thm-8}. Like before, $S$ must be a simple symmetric RW, but
now we also posit that the tail of the distribution of $\zeta_1$ is
\emph{regularly varying} with index $-\bq \le -1$. This means that
\begin{equation} \label{tail-zeta}
  \tau_\zeta(x) := \P( \zeta_1 > x ) = \frac{ \ell(x) } { x^\bq } ,
\end{equation}
where $\ell$ is a \emph{slowly varying} function at $+\infty$, namely, for all 
$c>0$,
\begin{equation}
  \lim_{x \to +\infty} \frac{ \ell(cx) } { \ell(x) } = 1.
  \footnote{See \cite{bgt} for a treatise on regularly varying functions.}
\end{equation}

In order to describe the upcoming theorems in their full power, we need
more notation. For $0<r<1$, set
\begin{equation}
  f_\bq(r) := \sum_{j=0}^{\lceil (1-r)/2r \rceil -1} \left( \left( \frac{2j+2}{1+r} 
  \right)^\bq - \left( \frac{2j}{1-r} \right)^\bq \right).
\end{equation}
It can be seen \cite[\S2.1]{z} that $0 < f_\bq(r) \le r^{-\bq}$ and, as $r \to 0^+$,
\begin{equation}
  f_\bq(r) \sim \frac1 {(\bq+1) \, r^\bq }.
\end{equation}
Here and in the rest of the paper $\sim$ denotes \emph{exact} asymptotic 
equivalence. This limit shows in particular that
$\int_0^1 r^{q-1} f_\bq(r) \, dr$ converges for $q=2\bq-1$ and $\bq>1$, or
$q>2\bq-1$ and $\bq \ge 1$. 
\begin{theorem} \label{thm-9}
  Assume that $S$ is a simple symmetric RW, $\mu < \infty$ and 
  $\tau_\zeta$ is regularly varying with index $-\bq \le -1$, cf.\ (\ref{tail-zeta}).
  For all $q>0$, recall the notation $\mq$ for the $q$-th moment of the 
  standard Gaussian, cf.\ (\ref{mq}), and set
  \begin{displaymath}
    d_{\mu, \bq, q} := \sqrt{\frac 2 \mu} \, \frac{\Gamma(q-\bq+1)} 
    {\Gamma(q-\bq+3/2)} \int_0^1 r^{q-1} f_\bq(r) \, dr.
  \end{displaymath}
  Then, as $t \to \infty$,
   \begin{displaymath}
     \E{|X(t)|^q} \sim \left\{
     \begin{array}{lll}
       \mq \, \mu^{q/2} \, t^{q/2} \,, && q < 2\bq-1, \text{ or } q=\bq=1; \\[2pt]
       \mq \, \mu^{q/2} \, t^{q/2} + d_{\mu, \bq, q} \, t^{q+1/2} \, \tau_\zeta(t) \,, 
         && q = 2\bq-1, \text{ and } \bq>1; \\[2pt]
       d_{\mu, \bq, q} \, t^{q+1/2} \, \tau_\zeta(t) \,, && q > 2\bq-1.
     \end{array}
    \right.
  \end{displaymath}
\end{theorem}

A few words of comment: In the subject of anomalous diffusion, an important 
quantity to investigate is the \emph{scaling exponent} of the moments,
\begin{equation}
  \gamma(q) := \lim_{t \to \infty} \frac{\log \E{|X(t)|^q}} {\log t},
\end{equation}
assuming this limit exists at least for a.e.\ $q >0$. In many relevant models
one observes that $q \mapsto \gamma(q)$ is piecewise linear with two branches,
a left one with slope $1/2$ and a right one with slope $1$. Researchers named
this situation \emph{strong anomalous diffusion},%
\footnote{Though different authors use different terminologies, not always 
compatible with each other, or even fully self-consistent.}
cf.\ \cite{cmmv, krs}. Theorem \ref{thm-9} shows that this is precisely what 
happens for the annealed L\'evy-Lorentz gas, under the above assumptions.
The corner between the two branches occurs at the moment of order 
$2\bq-1 \ge 2$, so the behavior of the second moment is still normal, at least
in terms of the leading exponent. Even more interestingly, this picture is very 
different from that of the corresponding quenched L\'evy-Lorentz gas, which 
is fully diffusive, as seen in Theorem \ref{thm-7}.%
\footnote{This does not mean that Theorems \ref{thm-7} and \ref{thm-9} are 
incompatible: what is happening here is that $t^{-q/2} \, E_\o[ |X(t)|^q ]$ 
converges to the suitable limit for a.a.\ $\o$, but, at least for large $q$, the 
convergence rate depends heavily on $\o$. Moreover, the convergence is 
not monotonic in $t$. Mathematically speaking, the convergence is neither
dominated nor monotonic, so one cannot interchange the limit in $t$ and
the integration on $\o$, to obtain the limit of the annealed moments
from that of the quenched moments.}

As for the quenched case, Theorem \ref{thm-9} is based on a large deviation 
result, which, however, is very different from Theorem \ref{thm-8}:

\begin{theorem} \label{thm-10}
  Under the same assumptions as in Theorem \ref{thm-9}, for $a \in (0,1]$,
  let
  \begin{displaymath}
    F_{\mu, \bq, a} := \frac1 {\sqrt{2\pi \mu}} \int_a^1 f_\bq \!\left( \frac a \eta \right)
    \frac{\eta^{-\bq}} {\sqrt{1-\eta}} \, d\eta.
  \end{displaymath}
  Then, as $t \to \infty$,
  \begin{displaymath}
    \P( X(t) > at ) = \P( X(t) < -at ) \sim F_{\mu, \bq, a} \sqrt{t} \, \tau_\zeta(t) .
  \end{displaymath}
  Moreover, for any $\delta \in (0,1)$, the lower order terms are uniformly 
  bounded for $a \in [\delta,1]$.
\end{theorem}

\begin{remark}
  The first equality of the above assertion is obvious because, by the 
  symmetry of the distributions of $\o$ and $S$, the annealed distribution
  of $X(t)$ is the same as that of $-X(t)$.
\end{remark}

\subsection{Infinite mean distance between targets, annealed theorems}

The last and hardest case is that of the CTRW $X$ in a medium $\o$ with 
$\mu=\infty$. At least to this author's knowledge, only a limit theorem 
seems to be available, that of Bianchi \emph{et al}, recently appeared in
\cite{blp}. We present it after some preparatory material.

First off, the assumption on the medium is that $\zeta_1$ is in the normal
domain of attraction of an $\a$-stable positive variable, with $\a \in (0,1)$.
As far $\o$ is 
concerned, this is the same assumption as in \S\ref{subs-y-inf-a}, so we 
use the same notation introduced there, in particular for the processes
$(Z^{(\a)}_{\pm}(s), \, s \ge 0)$ and $(Z^{(\a)}(s), \, s \in \R)$. The underlying 
random walk $S$ is assumed to be centered and such that $\E{|\xi_1|^q}
< \infty$, for some $q > 2/\a$. This implies in particular that $V_\xi = 
\E{|\xi_1|^2} < \infty$. So this is a special case of the assumptions on $S$
of \S\ref{subs-y-inf-a} (which were the same as in \S\ref{subs-y-fin-a}). 

All the preliminary results seen earlier then 
apply, in particular, for $n \to \infty$, $\diffusionscaledenvironmentnomean 
\dconv \diffusionlimitenvironment$, in $(\cadlag, J_1)$, cf.\ (\ref{hat-o-n}),
and $\diffusionscaledRWnomean := n^{-1/2} \, S_{\lfloor n \cdot\rfloor} \dconv 
W^{(2)}$, in $(\cadlag^+, J_1)$. Here $W^{(2)}$ is a Brownian motion 
such that $W^{(2)}(t)$ has mean $0$ and variance $V_\xi t$.
As clarified in \S\ref{subs-y-inf-a}, the processes $Z^{(\a)}_{\pm}$ and
$W^{(2)}$ are independent. Recalling the notation $M_\xi := \E{|\xi_1|}$, let 
$\Delta := (\Delta(t), t \ge 0)$ be defined by
\begin{equation} 
  \Delta(t) :=  M_\xi \left(\int_0^{\infty} L_t(x) \, d Z^{(\a)}_+ \! (x) + 
  \int_0^{\infty} L_t(-x) \, d Z^{(\a)}_- \! (x) \right) ,
\end{equation} 
where, for all $x\in\R$, $L_t(x) := \# \big( W^{(2)} |_{[0,t]} \big)^{-1} (x)$. In 
other words, $L_t(x)$ is the local time of the Brownian motion $W^{(2)}$ in 
$x$, up to time $t$. As a function of $x$, $L_t$ is compactly supported and 
almost surely continuous, thus the above r.h.s.\ is well-defined. Since
$L_t$ is also strictly increasing in $t$, $\Delta$ is almost surely 
continuous and strictly increasing. Processes 
like $\Delta$ are called \emph{Kesten-Spitzer processes} and 
arise in the context of \emph{RW in random scenery} \cite{ks}, which is
one of the technical ingredients of Theorem 2.1 of \cite{blp} which
we now present.
\begin{theorem} \label{thm-11}
  Under the above assumptions, in particular $\a \in (0,1)$, $\nu=0$, and 
  $\xi_1$ has a finite absolute moment of order $q > 2/\a$, let 
  $\hat{X}^{(n)} (t) := \ds \frac{X(nt)} {n^{1/(\a+1)}}$, for $t\ge 0$. Then the 
  annealed finite-dimensional distributions of $\hat{X}^{(n)}$ converge 
  to those of $Z^{(\a)} \circ W^{(2)} \circ \Delta^{-1}$. This means that, for all 
  $m \in \Z^+$ and $t_1,\ldots, t_m \in \R^+$,
  \begin{displaymath} 
    \big( \hat{X}^{(n)} (t_1), \ldots, \hat{X}^{(n)} (t_m) \big) \dconv \big( 
    Z^{(\a)} ( W^{(2)} ( \Delta^{-1}(t_1))), \ldots, Z^{(\a)} ( W^{(2)} ( 
    \Delta^{-1}(t_m))) \big) ,
  \end{displaymath} 
  as $n \to\infty$, relative to $\P$.
\end{theorem}

\begin{remark}
  The process $Z^{(\a)} \circ W^{(2)} \circ \Delta^{-1}$ is not a.s.\ c\`adl\`ag,
  so the same considerations and reference as in Remark \ref{rk-d} 
  apply here.
\end{remark}

\section{A brief discussion on perspectives}

Just by looking at the titles of the previous subsections, one notices that 
no quenched theorems were given  for the case of infinite mean 
distance between targets. This is the main shortcoming of the current
mathematical description of the processes $X$ and $Y$. Technically 
speaking, the problem is that, without the condition $\mu < \infty$, one
does not have a strong law of large numbers for the variables $\zeta_k$.
This is after all the simplest form of a quenched result and provides the 
scaling of $k \mapsto \o_k$, as $|k| \to \infty$, for \emph{each} realization 
$\o$ of the medium, apart from a negligible set of exceptions. How to prove 
quenched limit theorems without this basic ingredient is not clear to me at the 
moment.

An open question of a different nature is that of devising a 
good model of L\'evy-Lorentz gas in dimension $d\ge 2$. Here `good' 
means that it should have the following features, in one form or another:
\begin{itemize}
\item The random medium should be \emph{homogeneous}, in the sense 
that the distribution of the relative positions of two or more targets should not 
depend on the absolute position of any of the targets involved.%
\footnote{The reader who feels this condition is not well-defined is right, see 
footnote below.}

\item The distances between targets should be heavy-tailed.%
\footnote{This condition is ill-defined in the same way as the previous condition 
was. It would be well-defined if the points of the random medium were 
labeled in a consistent way, so that it would make sense to consider, say, the
distribution of the distance between $\o_k$ and $\o_\ell$ (here $k$ and $\ell$ 
are generic indices, not necessarily in $\Z$). But no labeling is assumed on the 
random medium, as it is not easy to think of a general, physically relevant way 
to label the points of a $d$-dimensional point process, for $d \ge 2$.}

\item The law of the random medium should be rotation-invariant, at least 
for a subgroup of rotations, e.g., the coordinate directions. In other words,
the model should be \emph{isotropic}, unless it has a clear reason not to be.

\item The transition probabilities from one target to the next should not 
depend on their distance, only on the ``degree of accessibility'' of the new 
target. For example, the next target might always be the nearest one along a 
random (isotropically selected) direction. The meaning of this condition is that it 
should be the medium, not the walker, to decide how long the next inertial 
stretch will be.
\end{itemize}

Even with all these features, a model might not be very interesting. Here is an 
example of a feasible, yet not very instructive model. Let 
$(\o'_{k_1}, \, k_1 \in \Z)$ and $(\o''_{k_2}, \, k_2 \in \Z)$ be two i.i.d.\ point 
processes in $\R$, as introduced in Section \ref{sec-intro}. For 
$k = (k_1, k_2) \in \Z^2$, set $\o_k := (\o'_{k_1}, \o''_{k_2})$. This defines a 
random medium $\o = (\o_k, \, k \in \Z^2)$ in $\R^2$.  An independent, 
\emph{$\Z^2$-valued}, underlying RW $(S_n, \, n \in \N)$ is  given, whereby we 
introduce the DTRW $Y := (Y_n := \o_{S_n}, \, n \in \N)$.  The CTRW 
$X := ( X(t), \, t \ge 0)$ is then defined as the unit-speed interpolation of $Y$.

Obviously, the process $X$ is simply the direct sum, in a very natural sense, 
of two independent, orthogonal, 1D CTRWs $X'$ and $X''$. Its properties are 
thus (for the most part) easily derived from those of $X'$ and $X''$, as 
presented in Section \ref{sec-ctrw}.

Introducing and investigating more relevant, and truly $d$-dimensional, flights
and walks in L\'evy random medium will be the subject of future work.

\section*{Acknowledgements}

I am indebted to all my coauthors in publications \cite{bcll, blp, sbblm}:
G.~Bet, A.~Bianchi, G.~Cristadoro, M.~Ligab\`o, E.~Magnanini, F.~P\`ene 
and S.~Stivanello. I acknowledge partial support by the PRIN Grant 
2017S35EHN \emph{``Regular and stochastic behaviour in dynamical 
systems''}  (MUR, Italy). This work is also part of my activity within the Gruppo 
Nazionale di Fisica Matematica (INdAM, Italy).

My research on RWs in L\'evy random media is my contribution to the joint 
UniBo-UniFi-UniPd project \emph{``Stochastic dynamics in disordered media 
and applications in the sciences''}. While the paper was in preparation, the 
Scientific Coordinator of the project in Florence, my dear friend and excellent 
mathematician Francesca Romana Nardi, passed away, leaving in all of us who 
were lucky enough to know her a vast sense of emptiness. This paper, while 
honoring the legacy of Carlo Cercignani, is dedicated to her memory.

\end{document}